\documentclass[11pt]{article}
\usepackage{amsmath,amssymb,theorem,amstext,amsgen,amsbsy,amsopn,amsfonts,graphicx,cases}

\usepackage{graphicx}
\usepackage{subfigure}
\usepackage{psfrag}
\usepackage{color}
\usepackage{enumerate}
\usepackage{epstopdf}
\usepackage[numbers,sort&compress]{natbib}
\usepackage{array}
\usepackage{ulem}
\usepackage{latexsym}
\newtheorem{thm}{Theorem}[section]
\newtheorem{conj}[thm]{Conjecture}
\newtheorem{cor}[thm]{Corollary}
\newtheorem{lem}[thm]{Lemma}
\newtheorem{prop}[thm]{Proposition}
\theorembodyfont{\rmfamily}

\def\pf{\bigskip\noindent {\bf Proof.}}

\def\dfn#1{{\sl #1}}

\def\less{\setminus}

\def\pf{\bigskip\noindent {\emph{Proof.}}}

\def\qed{ \hfill\vrule height3pt width6pt depth2pt}

\def\pf{\bigskip\noindent {\bf Proof.  }}

\title{Planar Tur\'an numbers of cubic graphs and disjoint union of cycles}

\author{Yongxin Lan\thanks{School of Science, Hebei University of Technology, Tianjin 300401,   China. Supported by the National Natural Science Foundation of China under Grant Nos. 12001154, 12161141006, Natural Science Foundation of Hebei Province under Grant No. A2021202025.   E-mail: yxlan@hebut.edu.cn}\hskip 1cm Yongtang Shi\thanks{Center for Combinatorics and LPMC, Nankai University, Tianjin  300071, China. Supported by the National Natural Science Foundation of China under Grant Nos. 11922112, 12161141006, Natural Science Foundation of Tianjin under Grant Nos. 20JCZDJC00840, 20JCJQJC00090.  E-mail: shi@nankai.edu.cn} \hskip 1cm
Zi-Xia Song\thanks{Department  of Mathematics, University of Central Florida, Orlando, FL 32816, USA.  Supported by the National   Science  Foundation under Grant No. DMS-1854903. E-mail:  Zixia.Song@ucf.edu}
}

\date{}

\begin{document}
\maketitle
\begin{abstract}
The planar Tur\'an number of  a graph  $H$, denoted  $ex_{_\mathcal{P}}(n,H)$,  is the maximum number of edges in a planar graph on $n$ vertices without containing   $H$ as a subgraph.   This  notion was introduced  by Dowden in 2016 and  has  attracted quite some attention since then;   those work mainly  focus on finding $ex_{_\mathcal{P}}(n,H)$ when $H$ is  a  cycle or  Theta graph or $H$ has maximum degree at least four.  In this paper, we  study    $ex_{_\mathcal{P}}(n,H)$ when $H$ is a cubic graph or    disjoint union of cycles or $H=K_{s, t}$.

\end{abstract}



\baselineskip 16pt
\section{Introduction}

All graphs considered in this paper are finite and simple.  We use $K_n$,    $C_n$  and $P_n$ to denote the complete graph, cycle and path on $n$ vertices, respectively.  Given a graph $G$, we  use  $|G| $ to denote  the number of vertices, $e(G)$ the number of edges, $\delta(G)$ the minimum degree, $\Delta(G)$ the maximum degree.     For a vertex $v\in V(G)$, we  define $N_G[v]:=N_G(v)\cup \{v\}$. For any $S,S'\subseteq V(G)$, we use $e_G(S,S')$ to denote the size of edge set $\{xy\in E(G)\mid x\in S \text{ and } y\in S'\}$.
 For any set  $S \subset V(G)$,
the subgraph of $G$ induced on  $S$, denoted $G[S]$, is the graph with vertex set $S$ and edge set $\{xy \in E(G) \mid  x, y \in S\}$. We denote by   $G \less S$ the subgraph of $G$ induced  on
$V(G) \less S$.   If $S=\{v\}$, then we simply write   $G\less v$.
The {\dfn{join}} $G+H$ (resp. {\dfn{union}} $G\cup H$)
of two vertex-disjoint graphs
$G$ and $H$ is the graph having vertex set $V(G)\cup V(H)$  and edge set $E(G)
\cup E(H)\cup \{xy\, |\,  x\in V(G),  y\in V(H)\}$  (resp. $E(G)\cup E(H)$). For a positive integer $t$ and a graph $H$, we use $tH$ to denote the disjoint union of $t$ copies of   $H$. Let $T_n$ denote a plane triangulation on $n\ge3$ vertices. We use  $T_n^-$ and $K_n^-$ to denote the graph   obtained from $T_n$  and $K_n$ with one edge removed, respectively. Given two isomorphic graphs $G$ and $H$,   we may (with a slight but common abuse of notation) write $G = H$.  For any positive integer $k$,  we define   $[k]:=\{1,2, \ldots, k\}$.

   \medskip

Given a graph $H$,  a graph is $H$-free if it does not contain   $H$ as a subgraph.   One of the best known results in extremal graph theory is Tur\'an's Theorem~\cite{Turan1941}, which
gives the maximum number of edges that a $K_t$-free graph on $n$ vertices can have. The celebrated  Erd\H{o}s-Stone Theorem~\cite{ErdosStone}  then extends this to the case when $K_t$ is
replaced by an arbitrary graph $H$ with at least one edge, showing that the maximum number of edges possible
is $(1+o(1)){n\choose 2}\left(\frac{\chi(H)-2}{\chi(H)-1}\right)$, where $\chi(H)$ denotes the chromatic number of $H$.  \medskip

In this paper, we continue to study the topic of ``extremal" planar graphs, that is,   how many edges
can an $H$-free planar graph on $n$ vertices have?   We define  $ex_{_\mathcal{P}}(n,H)$ to be  the maximum number of edges in an $H$-free planar graph on $n $ vertices.  Dowden~\cite{Dowden2016} initiated the study of $ex_{_\mathcal{P}}(n,H)$ and proved the following result.

 \begin{thm}[Dowden~\cite{Dowden2016}]\label{Dowden2016}  Let $n$ be a positive integer.
\begin{enumerate}[(a)]
\item $ex_{_\mathcal{P}}(n, C_3)=2n-4$ for all $n\geq 3$.
\item $ex_{_\mathcal{P}}(n, K_4)=3n-6$ for all $n\geq 4$.
\item $ex_{_\mathcal{P}}(n, C_4)\leq {15}(n-2)/7$ for all $n\geq 4$, with equality when $n\equiv 30 (\rm{mod}\, 70)$.
\item  $ex_{_\mathcal{P}}(n, C_5)\leq 12(n-2)/{5}$ for all $n\geq 5$.
\item  $ex_{_\mathcal{P}}(n, C_5)\leq (12n-33)/{5}$ for all $n\geq 11$. Equality holds for infinity many $n$.
\end{enumerate}
\end{thm}

 This topic has    attracted quite some attention since then. We refer the reader to a recent survey~\cite{LSS21} of the present authors for more information. Let $\Theta_k$ denote the family of Theta graphs on $k\ge4$ vertices, that is,   graphs obtained  from $C_k$    by  adding  an additional edge joining two non-consecutive vertices.    The present authors~\cite{LSS}   obtained tight upper bounds for  $ex_{_\mathcal{P}}(n, \Theta_k)$ for $k\in\{4,5\}$ and an upper bound for $ex_{_\mathcal{P}}(n, \Theta_6)$.
 \begin{thm}[Lan, Shi and Song~\cite{LSS}]\label{Theta}  Let $n$ be a positive integer.
\begin{enumerate}[(a)]
 \item $ex_{_\mathcal{P}}(n, \Theta_4)\le {12(n-2)}/5$ for all $n\geq 4$, with  equality when $n\equiv 12 (\rm{mod}\, 20)$.
 \item  $ex_{_\mathcal{P}}(n, \Theta_5)\le {5(n-2)}/2$ for all $n\ge5$,  with  equality when $n\equiv 50 (\rm{mod}\, 120)$.
\item $ex_{_\mathcal{P}}(n, C_6)\leq ex_{_\mathcal{P}}(n, \Theta_6) \le {18(n-2)}/7$ for all $n\geq6$.
\end{enumerate}
\end{thm}

Theorem~\ref{Theta}(c) has been strengthened by the authors in \cite{GGMPX20, Theta6} with tight upper bounds.
 \begin{thm}[Ghosh et al.~\cite{GGMPX20, Theta6}]\label{C6Theta6}  Let $n$ be a positive integer.
\begin{enumerate}[(a)]

\item $ex_{_\mathcal{P}}(n, C_6)\leq (5n-14)/2$ for all $n\geq 18$, with  equality when $n\equiv 2 (\rm{mod}\, 5)$.

\item $ex_{_\mathcal{P}}(n, \Theta_6)\leq (18n-48)/7$ for all $n\geq 14$. Equality  holds for infinitely many $n$.

\end{enumerate}
\end{thm}

As observed in \cite{Dowden2016}, for all  $n\ge6$,
 the planar triangulation $2K_1+C_{n-2}$ is $K_4$-free.  Hence,  $ex_{_\mathcal{P}}(n,H)=3n-6$ for all graphs $H$ which contains $K_4$ as a subgraph and   $n\ge \max\{|H|, 6\}$. The present authors~\cite{19LSS} also  investigated   a variety of sufficient conditions on $K_4$-free planar graphs $H$ such  that $ex_{_\mathcal{P}}(n,H)=3n-6$ for all $n\ge|H|$.

\begin{thm}[Lan, Shi and Song~\cite{19LSS}]\label{3n-6}
Let $H$ be a $K_4$-free planar graph and let  $n\ge  |H|$ be an integer. Then $ex_{_\mathcal{P}}(n,H)=3n-6$ if one of the following holds, where $n_{_k}(H)$  denotes  the number of  vertices   of degree $k$ in   $H$ for a positive integer $k $.

\begin{enumerate}[(a)]
\item   $\chi(H)=4$ and $n\ge|H|+2$.\label{g}
\item $\Delta(H)\ge 7$.\label{a}

\item $\Delta(H)=6$ and either $n_{_6}(H)+n_{_5}(H)\ge 2$ or  $n_{_6}(H)+n_{_5}(H)=1$ and
      $n_{_4}(H)\ge5$.\label{c}
\item $\Delta(H)=5$ and  either $H$ has   at least three $5$-vertices or $H$ has exactly two adjacent $5$-vertices. \label{d}
\item $\Delta(H)=4$ and  $n_{_4}(H)\ge7$. \label{f}

\item $H$ is $3$-regular  with $|H|\ge9$ or $H$ has at least three vertex-disjoint cycles or $H$ has exactly one vertex $u$ of degree $\Delta(H) \in\{4,5,6\}$  such that   $\Delta(H[N_H(u)])\ge 3$.\label{e}
\item   $\delta(H)\ge 4$ or   $H$ has exactly  one  vertex of degree at most $3$.\label{b}

\end{enumerate}

\end{thm}

In the same paper,  the present authors~\cite{19LSS} also   determined the values of $ex_{_\mathcal{P}}(n, H)$ when $H$ is a star, or wheel or $(t,r)$-fan.  Ghosh,  G\H{y}ori, Paulos and  Xiao~\cite{doublestar} recently determined the values of $ex_{_\mathcal{P}}(n,H)$ when $H$ is a double star.
Theorem~\ref{3n-6} implies that  $ex_{_\mathcal{P}}(n,H)$ remains wide open when $H$ is subcubic.  In particular, it seems quite non-trivial to determine $ex_{_\mathcal{P}}(n, C_k)$ for all $k\ge7$. Very recently, Cranston, Lidick\'y, Liu and Shantanam~\cite{CLLS21} proved that for each $k\ge 11$ and $n$ sufficiently large,
\[ex_{_\mathcal{P}}(n, C_k) >\left(3-\frac3k\right)n-6-\frac6k.\]
They further proposed the following conjecture.

 \begin{conj}[Cranston, Lidick\'y, Liu and Shantanam~\cite{CLLS21}]
There exists a constant $D$ such that for all $k$ and for all sufficiently large $n$, we have
$$ex_{_{\mathcal{P}}}(n,C_k)\le \left(3-\frac{3}{Dk^{lg_2^3}}\right)n.$$
\end{conj}

We propose a weaker conjecture here for the upper bound  of $ex_{_\mathcal{P}}(n,C_k)$ for all $k\ge 3$.
\begin{conj}\label{weak}
Let $n,k$ be integers with $n\ge k\ge3$. Then
$$ex_{_{\mathcal{P}}}(n,C_k)\le \left(3-\frac1{k-2}\right)n-4.$$
\end{conj}

   Conjecture~\ref{weak} holds for $3\le k\le 6$ and $n\ge k$ by Theorem~\ref{Dowden2016}(a,c,d) and Theorem~\ref{Theta}(c).  The assumed truth of Conjecture~\ref{weak}  will enable us  to prove in Section~\ref{disjointCk} that  $ex_{_{\mathcal{P}}}(n,tC_k\cup C_k^+)=ex_{_{\mathcal{P}}}(n,(t+1)C_k)$ for all $t\ge1$, $k\ge3$ and $n\ge (t+1)k+1$, where  $C_k^+$ denotes the graph on $k+1$ vertices obtained from $C_k$ by adding one pendant edge.   \medskip

The remainder of the  paper is organized as follows: we determine the values of $ex_{_\mathcal{P}}(n, H)$ when $H$ is $k$-regular in Section~\ref{cubic}; $H$ is disjoint union of cycles in Section~\ref{disjointCk}; $H=K_{s, t}$ in Section~\ref{Kst}. We end the paper with some concluding remarks in Section 5, including an improved lower bound for $ex_{_{\mathcal{P}}}(n,2C_k)$ for all $k\ge 7$.

 \section{Planar Tur\'an number of regular graphs}\label{cubic}

We begin this section with a lemma that will be essential in determining the planar Tur\'an numbers of regular planar graphs.

\begin{figure}[htbp]
\centering
\includegraphics*[scale=0.35]{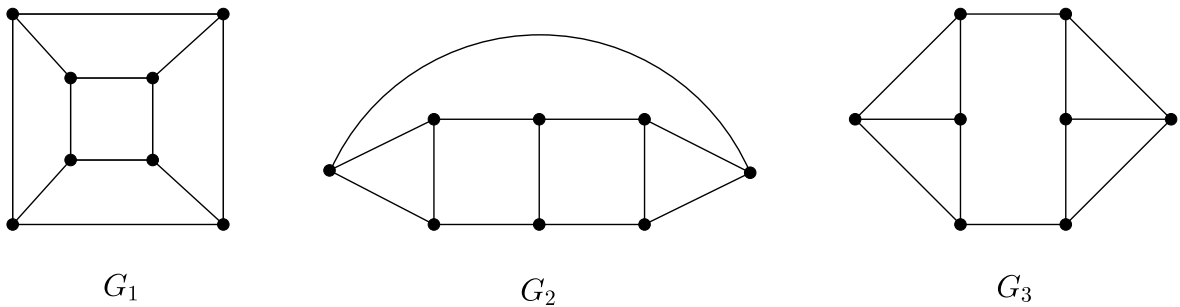}
\caption{Graphs $G_1$, $G_2$ and $ G_3$ }\label{G123}
\end{figure}

\begin{lem}\label{88}
Let  $G$ be a cubic planar graph on $8$ vertices. If $G$ is $K_4$-free, then $G\in\{G_1,G_2,G_3\}$, where   graphs $G_1, G_2, G_3$  are depicted in Figure~\ref{G123}.
\end{lem}

\pf
Let $G$ be given as in the statement. To establish the desired result, we first assume that there exists a vertex, say $v_1\in V(G)$, such that $v_1$ belongs to no triangle in $G$. Then $e(G[N_G[v_1]])=3$. Let $N_G(v_1)=\{v_2,v_3,v_4\}$ and $A=V(G)\less N_G[v_1]$.  Then $e_G(N_G(v_1),A)=6$ because $G$ is a cubic graph. Hence, $e(G[A])=e(G)-e(G[N_G[v_1]])-e_G(N_G(v_1),A)=3$,  which implies that either $G[A]= K_{1,3}$ or $G[A]= P_{4}$. If $G[A]= K_{1,3}$ and $x\in A$ is the center of $K_{1,3}$, then $G\less \{v_1,x\} = C_6$ and so $G=G_1$. So we next assume that $G[A]= P_{4}$. Let $G[A]$ be a path with vertices $v_5,v_6,v_7,v_8$ in order. We claim that $v_6$ and $v_7$ have no common neighbour in $G$. Suppose $v_6$ and $v_7$ have a common neighbour, say $v_2$. Then $N_G(v_3)=N_G(v_4)=\{v_1,v_5,v_8\}$ because $G$ is a cubic graph. But then $G$ has a $K_{3,3}$-minor with one part $\{v_2,v_3,v_4\}$ and the other part $\{v_1,\{v_5,v_6\},\{v_7,v_8\}\}$, a contradiction. Thus, $v_6$ and $v_7$ have no common neighbour in $G$. Without loss of generality, we may assume that $v_6v_2,v_7v_3\in E(G)$. Notice that $N_G(v_4)=\{v_1,v_5,v_8\}$ because $G$ is a cubic graph. Then we see that either $v_5v_2\in E(G)$ or $v_5v_3\in E(G)$. If $v_5v_3\in E(G)$, then $v_8v_2\in E(G)$ and so $G$ contains a $K_{3,3}$-minor with one part $\{v_2,v_3,v_4\}$ and the other part $\{v_1,\{v_5,v_6\},\{v_7,v_8\}\}$, a contradiction. Hence, $v_5v_2\in E(G)$ and so $v_8v_3\in E(G)$, which implies that $G=G_2$.

Next we assume that every vertex in $G$ belongs to at least one triangle. We claim that there must exist a vertex such that it belongs to two triangles. Suppose that every vertex in $G$ belongs to exactly one triangle. Let $v_1\in V(G)$. Then $G[N_G[v_1]]= K_1+(K_2\cup K_1)$. Let $A=V(G)\less N_G[v_1]$. Then $e_G(N_G(v_1),A)=4$ because $G$ is a cubic graph. Hence, $e(G[A])=e(G)-e(G[N_G[v_1]])-e_G(N_G(v_1),A)=4$, which implies that either $G[A]= C_4$ or $G[A]= K_1+(K_2\cup K_1)$. If $G[A]= C_4$, there exist two vertices in $A$ which does not belong to any triangle, a contradiction. If $G[A]= K_1+(K_2\cup K_1)$, then there exists one vertex such that it belongs to either two triangles or no triangle, a contradiction. Thus, there must exist a vertex belonging to two triangles. Assume that $v_1$ belongs to two triangles in $G$. Since $G$ is $K_4$-free, we have $G[N_G[v_1]]= K_4^-$. Let $A=V(G)\less N_G[v_1]$. Then $e_G(N_G(v_1),A)=2$ because $G$ is a cubic graph. Hence, $e(G[A])=e(G)-e(G[N_G[v_1]])-e_G(N_G(v_1),A)=5$, which implies that $G[A]= K_4^-$. Thus, $G=G_3$. \qed

\begin{thm}\label{main}
Let $H$ be a $k$-regular planar graph with $k\ge3$ and let  $n\ge  |H|$ be an integer. Then \[ex_{_\mathcal{P}}(n,H)=
\begin{cases}
\, 3n-6 & \text{ if } \, |H|\ge8, \text{ or } |H|=6 \text{ and } n\ge 10;\\[2mm]

\, 3n-7 & \text{ if } \, |H|=6 \text{ and } n\le 9.
\end{cases}
\]
\end{thm}
\pf
Let $H$ and $n$  be given as in the statement.  By Theorem~\ref{Dowden2016}($b$) and Theorem~\ref{3n-6}($f, g$),  we see that $ex_{_\mathcal{P}}(n,H)=3n-6$ for all $n\ge |H|$  if  $H$ contains a copy of $K_4$, or $k=3$ and $|H|\ge 9$,   or $k\ge4$. For the remainder of the proof, let $H$  be a  $K_4$-free cubic planar graph with $|H|\le8$.  Assume first that $|H|=8$. By Lemma~\ref{88}, we see that $H\in\{G_1,G_2,G_3\}$. Then the planar triangulation $2K_1+C_{n-2}$ is $H$-free when $n>|H|$, and the planar triangulation $K_2+P_{n-2}$ is $H$-free when $n=|H|$. Hence, $ex_{_\mathcal{P}}(n,H)=3n-6$ for all $n\ge |H|$ and $|H|=8$.\medskip

\begin{figure}[htbp]
\centering
\includegraphics*[scale=0.4]{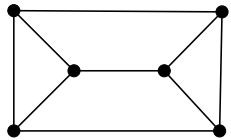}
\caption{The unique cubic planar graph $H$ on $6$ vertices.}\label{6}
\end{figure}

\begin{figure}[htbp]
\centering
\includegraphics*[scale=0.4]{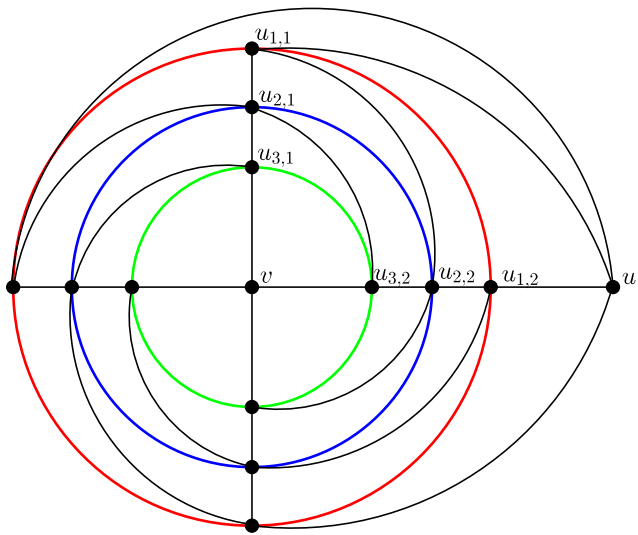}
\caption{The plane triangulation $Q_k$ when $k=3$, where $C^1,C^2,C^3$ are in red, blue and green, respectively.}\label{qk}
\end{figure}

It remains to consider the case when $|H|=6$. Then $H$ is the graph given in Figure~\ref{6}.
We next show that  $ex_{_\mathcal{P}}(n,H)=3n-6$ for all $n\ge 10$.  Let  $n:=4k+2+\ell$ for some $\ell\in\{0,1,2,3\}$ and integer  $k\ge2$.    Let $Q_k$ be a plane triangulation on $n=4k+2$ vertices constructed as follows: for each $i\in[k]$, let $C^i$ be a cycle with vertices $u_{i,1},u_{i,2},u_{i,3},u_{i,4}$ in order, let $Q_k$ be the plane triangulation obtained from disjoint union of $C^1,\ldots,C^k$ by adding edges $u_{i,j}u_{i+1,j}$ and $u_{i,j}u_{i+1,j+1}$ for all $i\in[k-1]$ and $j\in[4]$, where all arithmetic on the index $j+1$ here is done modulo 4, and finally adding two new nonadjacent vertices $u$ and $v$ such that $u$ is adjacent to all vertices of $C^1$ and $v$ is adjacent to all vertices of $C^k$. The graph $Q_k$ when $k=3$ is depicted in Figure~\ref{qk}.  Let $Q_k^\ell=Q_k$ if $\ell=0$. For $\ell\in\{1,2,3\}$, let $F_j$ be the face of $Q_k$ with vertices $u_{k-1,j},u_{k,j},u_{k,j+1}$ for each $j\in[\ell]$,  and let $Q_k^\ell$ be the plane triangulation on $n$ vertices obtained from $Q_k$ by adding one new vertex, say $x_j$,  adjacent to the three vertices on the boundary of  $F_j$ for each $j\in[\ell]$.  It can be checked that  $Q_k^\ell$ is   $H$-free. Therefore, $ex_{_\mathcal{P}}(n,H)=3n-6$ for all $n\ge 10$. \medskip

We next show that $ex_{_\mathcal{P}}(n,H)=3n-7$ for each $n\in\{6,7,8,9\}$. To obtain the desired upper bound, it suffices to show that every plane triangulation $T$ on $n\in\{6,7,8,9\}$ vertices contains a copy of $H$.  Let $x\in V(T)$ be a vertex of maximum degree in $T$. Then $d_T(x)\ge 4$ and $T[N_T[x]]$ is a wheel on $d_T(x)+1$ vertices. It is simple to check that   $T$  contains a copy of $H$, and so $ex_{_\mathcal{P}}(n,H)\le 3n-7$ for each $n\in\{6,7,8,9\}$.
 On the other hand, for each $n\in\{6,7, 8\}$, the planar graph $K_2+(P_3\cup P_{n-5})$    is $H$-free with  $3n-7$ edges;   for $n=9$, the planar graph $Q_2\less u$  is $H$-free with   $3n-7$ edges. Hence, $ex_{_\mathcal{P}}(n,H)= 3n-7$ for each $n\in\{6,7,8,9\}$, as desired.
\qed

\section{Planar Tur\'an number of disconnected graphs}\label{disjointCk}

  Given a plane graph $G$ and an integer $i\ge 3$, an $i$-face in $G$ is a face of order $i$.  Let $f_i$ and $f(G)$ denote the number of $i$-faces and all faces in $G$, respectively.
In this section we study the planar Tur\'an number of disconnected graphs.
We first consider $tC_k$, the $t$ vertex-disjoint copies of $C_k$,  and give a tight bound for $ex_{_\mathcal{P}}(n,tC_3)$ for all $n\ge3t\ge3$. It is worth noting that  $ex_{_\mathcal{P}}(n,H)=3n-6$   if $H$ has   three vertex-disjoint cycles, due to  Theorem~\ref{3n-6}($f$).

\begin{thm}\label{tC3}
For   integers $n$ and $t$ with $n\ge 3t\ge 3$, we have
\[ex_{_\mathcal{P}}(n,tC_3)=
\begin{cases}
3n-6 &\, \text{ if } \, t\ge3;\\
\left\lceil\frac{5n}{2}\right\rceil-5 &\, \text{ if } \, t=2;\\
2n-4 &\, \text{ if } \, t=1.\\
\end{cases}\]
\end{thm}

\pf By Theorem~\ref{Dowden2016}($a$) and Theorem~\ref{3n-6}($f$),    $ex_{_\mathcal{P}}(n,tC_3)=2n-4$ if $t=1$, and $ex_{_\mathcal{P}}(n,tC_3)=3n-6$ if $t\ge 3$.  We may assume that $t=2$. Then $n\ge6$. We first show that $ex_{_\mathcal{P}}(n,2C_3)\ge \lceil5n/2\rceil-5$. Let $P$ be a path on $n-2$ vertices and $S$ be a maximum independent set of $P$ containing the two ends of $P$.
 Let $G$ be the planar graph on $n $ vertices obtained from $P$ by adding two new adjacent vertices $u$ and $v$ such that $u$ is joined to  every vertex in  $V(P)$  and   $v$ is joined to  every vertex in $S$.
  Then $G$ is    $2C_3$-free    with $|G|=n$  and  $e(G)=(n-3)+(n-1)+\lceil(n-2)/2\rceil =\lceil 5n/2\rceil-5$. Hence, $ex_{_\mathcal{P}}(n,2C_3)\ge e(G)=\lceil5n/2\rceil-5$.\medskip

\begin{figure}[htbp]
\centering
\includegraphics*[scale=0.3]{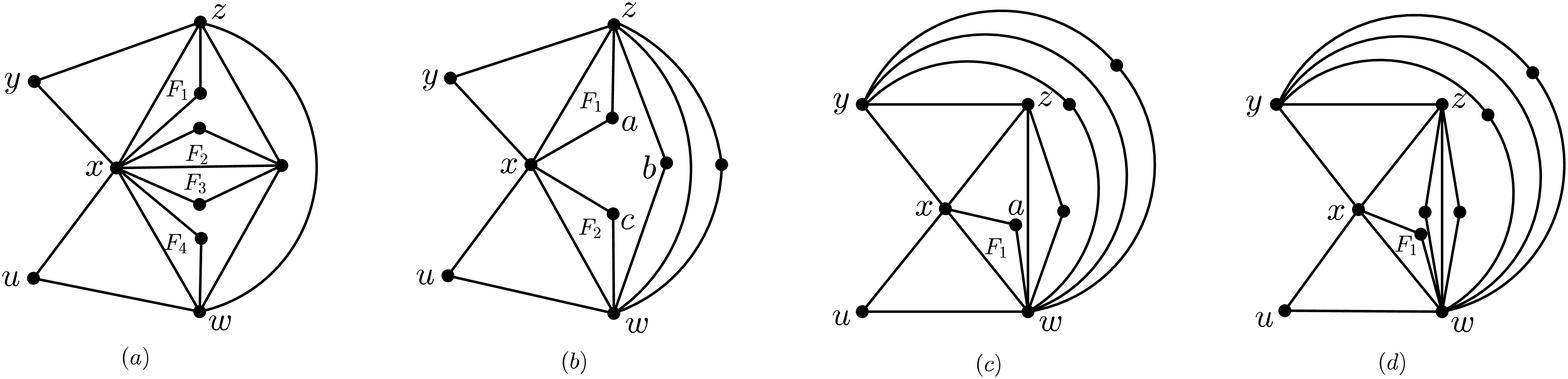}
\caption{$(a)$ $\mathcal{F}(x)=\{F',F'',F_1,F_2,F_3,F_4\}$; $(b)$ $\mathcal{F}(x)=\{F',F'',F_1,F_2\}$; $(c,d)$ $\mathcal{F}(x)=\{F',F'',F_1\}$.}\label{8}
\end{figure}

\begin{figure}[htbp]
\centering
\includegraphics*[scale=0.3]{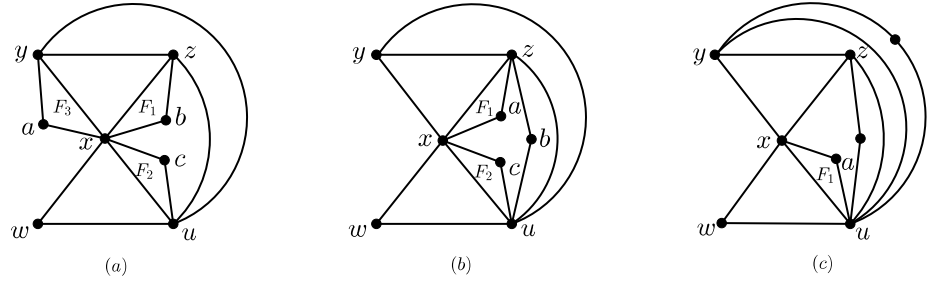}
\caption{$(a)$ $\mathcal{F}(x)=\{F',F'',F_1,F_2,F_3\}$; $(b)$ $\mathcal{F}(x)=\{F',F'',F_1,F_2\}$; $(c)$ $\mathcal{F}(x)=\{F',F'',F_1\}$.}\label{9}
\end{figure}

We next show that $ex_{_\mathcal{P}}(n,2C_3)\le \lceil5n/2\rceil-5$.  It can be easily checked that  every plane graph $T_n^-$ contains a copy of $2C_3$, where $n\in\{6,7\}$. Hence, $ex_{_\mathcal{P}}(n,2C_3)\le e(T_n^-)-1=3n-8=\lceil5n/2\rceil-5$ when $n\in\{6,7\}$. We may assume that $n\ge 8$. Let $G$ be a  $2C_3$-free plane graph on $n \ge8$ vertices.  We claim that $f_3\le  n-1$.
Suppose  $f_3\ge n\ge8 $.   Let $\mathcal{F}$ be the set of all $3$-faces of $G$. Then $|\mathcal{F}|=f_3$. For each $v\in V(G)$,  let  $\mathcal{F}(v):=\{F\in \mathcal{F}\mid v\in V(F)\}$. Then       $|\mathcal{F}(v)|\le n-1$ and so  $\mathcal{F}\setminus\mathcal{F}(v)\ne \emptyset$.
  Since $G$ is $2C_3$-free, we see that $V(F)\cap V(F')\ne \emptyset$ for every  pair $F, F'\in\mathcal{F}$.     Since $f_3\ge n\ge7$, there exist $F', F''\in \mathcal{F}$ such that $|V(F')\cap V(F'')|=1$. We may assume that  $V(F')=\{x,y,z\}$ and   $V(F'')=\{x,u,w\}$, where $x,y, z, u, w$ are pairwise distinct.     It follows that for every $F\in\mathcal{F}\setminus \mathcal{F}(x)$, we have $|V(F)\cap\{y,z\}|\ge1$ and $|V(F)\cap\{u,w\}|\ge1$.  Suppose  $|V(F)\cap\{y,z\}|=1$ and $|V(F)\cap\{u,w\}|=1$ for every $F\in\mathcal{F}\setminus \mathcal{F}(x)$.  Then $|\mathcal{F}\setminus \mathcal{F}(x)|\le4$. In addition, if $|\mathcal{F}\setminus \mathcal{F}(x)|=1$, then $|\mathcal{F}(x)|\le 6$ (see Figure~\ref{8}$(a)$ when $|\mathcal{F}(x)|=6$); if $|\mathcal{F}\setminus \mathcal{F}(x)|=2$, then $|\mathcal{F}(x)|\le 4$ (see Figure~\ref{8}($b$) when $|\mathcal{F}(x)|=4$, where vertices $a,b,c$ are not necessary distinct); if $3\le|\mathcal{F}\setminus \mathcal{F}(x)|\le4$, then $|\mathcal{F}(x)|\le 3$ (see Figure~\ref{8}$(c,d)$ when $|\mathcal{F}(x)|=3$, where vertices $a$ and $z$ in Figure~$(c)$ are not necessary distinct).
  Thus $f_3=|\mathcal{F}\setminus \mathcal{F}(x)| +|\mathcal{F}(x)|\le 7$,  contrary to the assumption that $f_3\ge n\ge 8$.
  Thus for some $F^*\in\mathcal{F}\setminus \mathcal{F}(x)$, we have    $y, z\in V(F^*)$ or $u, w\in V(F^*)$, say the former.
  We may further assume that $u\in V(F^*)$. Then $|\mathcal{F}\setminus \mathcal{F}(x)|\le3$. In addition, if $|\mathcal{F}\setminus \mathcal{F}(x)|=1$, then $|\mathcal{F}(x)|\le 5$ (see Figure~\ref{9}$(a)$ when $|\mathcal{F}(x)|=5$, where vertices $a$ and $w$, or $b$ and $c$ are not necessary distinct); if $2\le|\mathcal{F}\setminus \mathcal{F}(x)|\le3$, then $|\mathcal{F}(x)|\le 4$ (see Figure~\ref{9}$(b,c)$, where vertices $a,b,c$ in Figure~$(b)$ and $a,z$ in Figure~$(c)$ are not necessary distinct).
  It follows that $f_3=|\mathcal{F}\setminus \mathcal{F}(x)| +|\mathcal{F}(x)|\le 7$, contrary to the assumption that $f_3\ge n\ge8$.  This proves that $f_3\le  n-1$, as claimed.   \medskip


  Note that $G$ has no face of size at most two.  It follows that
$$2e(G)=3f_3+\sum_{i\ge4}if_i\ge3f_3+4(f(G)-f_3)=4f(G)-f_3\ge4f(G)-(n-1),$$
which implies that $f(G)\le (2e(G)+n-1)/4$. By Euler's formula, \[n-2=e(G)-f(G)\ge e(G)/2-(n-1)/4.\] Hence, $e(G)\le \lceil5n/2\rceil-5$, as desired. This completes the proof of Theorem~\ref{tC3}.\qed\bigskip

 We next investigate  $ex_{_\mathcal{P}}(n, tC_k\cup C_k^+)$  and $ex_{_\mathcal{P}}(n,(t+1)C_k)$ are the same, where $t\ge 1$ is an integer.
 We begin with a lower bound  construction for $ex_{_\mathcal{P}}(n,2C_k)$ which plays a key role in the proof of Theorem~\ref{tCk+}.

\begin{lem}\label{2Ck}
Let $n$ and $k$ be positive integers with  $n\ge 2k\ge8$. Let $r$ be the remainder of $n-3$ when divided by $k-2$.   Then
\[ex_{_\mathcal{P}}(n,2C_k)\ge\left(3-\frac{1}{k-2}\right)n+\frac{3+r}{k-2} -5+\max\{1-r, 0\}.\]
\end{lem}

\pf Let $n,k,r$ be given as in the statement. Let $t\ge0$ be an integer satisfying \[(2k-1)+t(k-2)+r=n.\] Let $P_1, P_2,\ldots,P_{t+1}$  be vertex-disjoint  paths with  $|P_{t+1}|=2k-3$ and $|P_i|=k-2$ for each $i\in [t]$.   Let $H=P_1\cup\cdots\cup P_{t+1}$. Then $|H|=(2k-3)+t(k-2)=n-2-r$ and \[e(H)=2k-4+t(k-3)=2k-4+(n-2k+1-r)\left(1-\frac1{k-2}\right).\] Assume first  $r=0$.  Let $G$ be  the planar  graph obtained from $H$ by adding two new adjacent vertices such that each new vertex is joined to each vertex of $H$. Clearly, $G$ is $2C_k$-free with $|G|=|H|+2=(n-2)+2=n$ and \[e(G)=2n-3+e(H)=2n-3+2k-4+(n-2k+1)\left(1-\frac1{k-2}\right)
=\left(3-\frac{1}{k-2}\right)n+\frac{3}{k-2}-4.\] Hence, $ex_{_\mathcal{P}}(n,2C_k)\ge e(G)=(3-\frac{1}{k-2})n+ \frac{3}{k-2}-4$.  We may assume that $r\ge1$. Let $Q$ be a path on $r$ vertices such that $V(Q)\cap V(H)=\emptyset$. Let $G'$ be the planar  graph obtained from $H\cup Q$ by adding two new adjacent vertices such that each new vertex is joined to each vertex of $H\cup Q$. Clearly, $G'$ is $2C_k$-free with $|G'|=|H|+| Q|+2=(n-2-r)+r+2=n$ and
\begin{align*}
e(G')&=2n-3+e(H)+e(Q)\\
&=2n-3+2k-4+(n-2k+1-r)\left(1-\frac1{k-2}\right)+(r-1)\\
&=\left(3-\frac{1}{k-2}\right)n+\frac{3+r}{k-2}-5.
\end{align*}
 Hence, $ex_{_\mathcal{P}}(n,2C_k)\ge e(G')=\left(3-\frac{1}{k-2}\right)n+\frac{3+r}{k-2}-5$, as desired. \medskip

 This completes the proof of Lemma~\ref{2Ck}.\qed\bigskip

   We are ready to prove that if Conjecture~\ref{weak} is true, then   $ex_{_\mathcal{P}}(n, tC_k\cup C_k^+)$  and $ex_{_\mathcal{P}}(n,(t+1)C_k)$ are the same  for all $n\ge (t+1)k+1$.

\begin{thm}\label{tCk+}
Let $t\ge1$, $k\ge3$ and $n\ge (t+1)k+1$ be integers.  If Conjecture~\ref{weak} is true, then
  \[ex_{_{\mathcal{P}}}(n,tC_k\cup C_k^+)=ex_{_{\mathcal{P}}}(n,(t+1)C_k).\]
\end{thm}

\pf Since    $(t+1)C_k$ is a subgraph of $tC_k\cup C_k^+$, we have $ex_{_{\mathcal{P}}}(n,tC_k\cup C_k^+)\ge ex_{_{\mathcal{P}}}(n,(t+1)C_k)$. We next  prove that $ex_{_{\mathcal{P}}}(n,tC_k\cup C_k^+) \le ex_{_{\mathcal{P}}}(n,(t+1)C_k)$.  By Theorem~\ref{3n-6}($f$), $ex_{_{\mathcal{P}}}(n,tC_k\cup C_k^+)=ex_{_{\mathcal{P}}}(n,(t+1)C_k)=3n-6$ when $t\ge2$. We may assume that $t=1$.
Let $G$ be any $C_k\cup C_k^+$-free planar graph with  $n:=|G|\ge 2k+1$ and  $e(G)=ex_{_{\mathcal{P}}}(n,C_k\cup C_k^+)$. If $G$ is $2C_k$-free, then $e(G)\le ex_{_{\mathcal{P}}}(n,2C_k)$, as desired.  We may assume that  $G$ contains   $2C_k$ as a subgraph.
Assume first that     $G$ contains a subgraph  $H^*:=2C_k$  such that $G[V(H^*)]$ is connected.   Then $G[V(H^*)]$ contains $C_k^+$ as a subgraph. Let $H:=G\less  V(H^*)$.  Then  $|H|\ge 1$ and $H$ is $C_k$-free.  Note that  $e(H)\le ex_{_{\mathcal{P}}}(|H|, C_k)\le   \left(3-\frac1{k-2}\right)|H|-4$ when      $|H|\ge k $ by the assumed truth of Conjecture~\ref{weak};   $e(H)\le \max\{3|H|-6, |H|-1\}$ when $1\le |H|\le k-1$.   Therefore,
\begin{align*}
e(G)=e(G[V(H^*)])+ e(H)&\le (3\times 2k-6)+e(H) \\
&\le \begin{cases}  \left(3-\frac1{k-2}\right)n-8+\frac4{k-2} &\mathrm{if} \,\,  |H|\ge k  \\
  3n-12   &\mathrm{if} \,\, 3\le |H|\le  k-1\\
   3n-7-2|H|   &\mathrm{if} \,\, 1\le |H|\le2.  \end{cases}
\end{align*}
 By Lemma~\ref{2Ck} when $k\ge4$ and Theorem~\ref{tC3} when $k=3$,    we have  $e(G)\le ex_{_{\mathcal{P}}}(n,2C_k)$, as desired.  We may assume that  $G$ contains no such subgraph $H^*$.
 Let $H_1,\ldots,H_{p} $ be  vertex-disjoint subgraphs of   $G$ such that $H_i:= C_k$ for each $i\in[p]$ and $H:=G\less \sum_{i=1}^p V(H_i)$ is $ C_k$-free, where $p\ge 2$ is an integer.
 Then   $|H|=n-kp $ and  no vertex in $V(G)\less V(H_i)$ is adjacent to any vertex in $V(H_i)$ for all $i\in[p]$.
 Then $e(G[V(H_i)])\le 3|H_i|-6=3k-6$ for all $i\in[p]$,  and    $e(H)\le ex_{_{\mathcal{P}}}(|H|, C_k)\le   \left(3-\frac1{k-2}\right)|H|-4$ when  $|H|\ge k$ by the assumed truth of Conjecture~\ref{weak};  $e(H)\le \max\{3|H|-6, |H|\}$ when $  |H|\le k-1$.   Note that    $|H|=n-kp\le n-2k $.
 Then
 \begin{align*}
e(G)&=e(G[V(H_1)])+\cdots+e(G[V(H_p)])+e(H)\\
 &= p(3k-6)+e(H)\\
&\le \begin{cases}  \left(3-\frac1{k-2}\right)n-4-p\left(5-\frac{2}{k-2}\right)\le  \left(3-\frac1{k-2}\right)n-14+\frac{4}{k-2} &\mathrm{if} \,\,  |H|\ge k  \\
  3n-18  &\mathrm{if} \,\, 3\le |H|\le  k-1\\
   3n-12-2|H|   &\mathrm{if} \,\,   |H|\le2.  \end{cases}
\end{align*}
By Lemma~\ref{2Ck} when $k\ge4$ and  Theorem~\ref{tC3}  when $k=3$,     we have $e(G)\le ex_{_{\mathcal{P}}}(n,2C_k)$, as desired.\medskip

This completes the proof of Theorem~\ref{tCk+}. \qed\bigskip

Corollary~\ref{c:tCk+} follows directly from Theorem~\ref{tCk+},
Theorem~\ref{Dowden2016}($a,c,d$) and Theorem~\ref{Theta}($c$).
\begin{cor}\label{c:tCk+}
Let $t\ge1$, $3\le k\le6$ and $n\ge (t+1)k+1$ be integers. Then
  \[ex_{_{\mathcal{P}}}(n,tC_k\cup C_k^+)=ex_{_{\mathcal{P}}}(n,(t+1)C_k).\]
\end{cor}

We end this section with a result showing that   $ex_{_\mathcal{P}}(n,   C_3^+)=ex_{_\mathcal{P}}(n,C_3)$   for all $n\ge 4$.
\begin{prop}\label{p:C3+}
For all $n\ge4$,  $ex_{_{\mathcal{P}}}(n,  C_3^+)=ex_{_{\mathcal{P}}}(n,C_3)$.
   \end{prop}

\pf Since  $C_3$ is a subgraph of $C_3^+$, we have $ ex_{_{\mathcal{P}}}(n,C_3^+)\ge ex_{_{\mathcal{P}}}(n,C_3)$. We next prove that $ex_{_{\mathcal{P}}}(n,C_3^+)\le ex_{_{\mathcal{P}}}(n,C_3)$.
Let $G$ be any $C_3^+$-free planar graph on $n\ge 4$ vertices with $e(G)=ex_{_{\mathcal{P}}}(n, C_3^+)$. If $G$ is $C_3$-free, then $e(G)\le ex_{_{\mathcal{P}}}(n,C_3)$, as desired.  We may assume that  $G$ contains   $C_3$ as a subgraph. Let $ H_1,\ldots,H_{p} $ be  vertex-disjoint subgraphs of   $G$ such that   $H_i=C_3$ for each $i \in[p]$ and $H:=G\less \sum_{i=1}^p V(H_i)$ is $C_3$-free, where $p\ge1$ is an integer.  Since    $G$ is  $  C_3^+$-free, it follows that $|H|=n-3p$ and  no vertex in $V(G)\less V(H_i)$ is adjacent to any vertex in $V(H_i)$ for all $i\in[p]$. Therefore, $e(G[V(H_i)])=3$ for all $i\in[p]$ and   $e(H)\le ex_{_{\mathcal{P}}}(|H|, C_3)=2|H|-4$ when  $|H|\ge 3$ by Theorem~\ref{Dowden2016}($a$);  $e(H)\le |H|$ when $  |H|\le 2$.   Note that    $|H|=n-3p\le n-3 $
 Then
 \begin{align*}
e(G)=e(H_1)+\cdots+e(H_p)+e(H)
\le 3p+e(H) =n-|H|+e(H)
\le \begin{cases}  2n-7 &\mathrm{if} \,\,  |H|\ge 3  \\
  n   &\mathrm{if} \,\, |H|\le2.  \end{cases}
\end{align*}
 By Theorem~\ref{Dowden2016}($a$),     $e(G)\le ex_{_{\mathcal{P}}}(n, C_3)$, as desired.\medskip

This completes the proof of Proposition~\ref{p:C3+}.\qed\medskip

\section{Planar Tur\'an number of complete bipartite graphs}\label{Kst}
Finally, we study   the planar Tur\'an number  of $K_{m,t}$, where $t\ge m\ge1$.  Note that  $ex_{_{\mathcal{P}}}(n, K_{m,t})=3n-6 $ when $  m\ge3$;  Theorem 6 in \cite{19LSS}  completely determines the values of $ex_{_{\mathcal{P}}}(n, K_{m,t})$ when  $m=1$; Theorem~\ref{Dowden2016}($c$) settles $ex_{_{\mathcal{P}}}(n, K_{m,t})$ for the case $m=t=2$.  We prove the remaining cases for $ex_{_{\mathcal{P}}}(n, K_{2,t})$.
Let $O_n$ denote the unique outerplane graph with  $2n-3$ edge, maximum degree 4, and the outer face of order $n$.
\begin{thm}\label{bipartite}
For     integers $t\ge3$ and $n\ge t+2$,
we have
\[ex_{_{\mathcal{P}}}(n, K_{2,t})=
\begin{cases}
3n-6 &\, \text{ if }\; t\ge5\; \text{ and }\;n\ge t+2, \;\text{ or } \;t=4 \;\text{ and } \;n\ge9,  \;\text{ or } \;t=3 \;\text{ and } \;n\ge12;\\
3n-7 &\, \text{ if }\; t=4 \; \text{ and }\; n\le 8;\\
3n-8 &\, \text{ if } \; t=3 \; \text{ and } \; n\le11.\\
\end{cases}
\]
\end{thm}
\vspace{-4mm}
\pf Let $n$ and $t$ be given as in the statement.
Note that $O_{n-1}$ is $K_{2,3}$-free. It follows that $K_1+O_{n-1}$ is $K_{2,5}$-free. Hence, $ex_{_{\mathcal{P}}}(n, K_{2,t})=3n-6$ for all $t\ge5$ and $n\ge t+2$.
We may  assume that $t\in\{3,4\}$. We first consider the case $n\ge12$. Let $H_{\ell}$ be the plane graph on $6\ell$ vertices constructed as follows: for each $i\in [\ell]$, let $C^i$ be a cycle with vertices $u_{i,1},u_{i,2},\ldots,u_{i,6}$ in order, let $H_{\ell}$ be the plane graph obtained from disjoint union of $C^1,\ldots,C^{\ell}$ by adding edges $u_{i,j}u_{i+1,j}$ and $u_{i,j}u_{i+1,j+1}$ for all $i\in[\ell-1]$ and $j\in[6]$, where all arithmetic on the index $j+1$ here   is done modulo 6. For $i\in\{1,\ell\}$, let
 $W_5^i$ be a wheel on $6$ vertices   $v_{i,0},v_{i,1},v_{i,2},v_{i,3},v_{i,4},v_{i,5}$, where $v_{i,0}$ is the center vertex of   $W_5^i$; $K_p^i$ be a complete graph with vertices $w_{i,1},w_{i,2},\ldots,w_{i,p}$.

\begin{figure}[htbp]
\centering
\includegraphics*[scale=0.25]{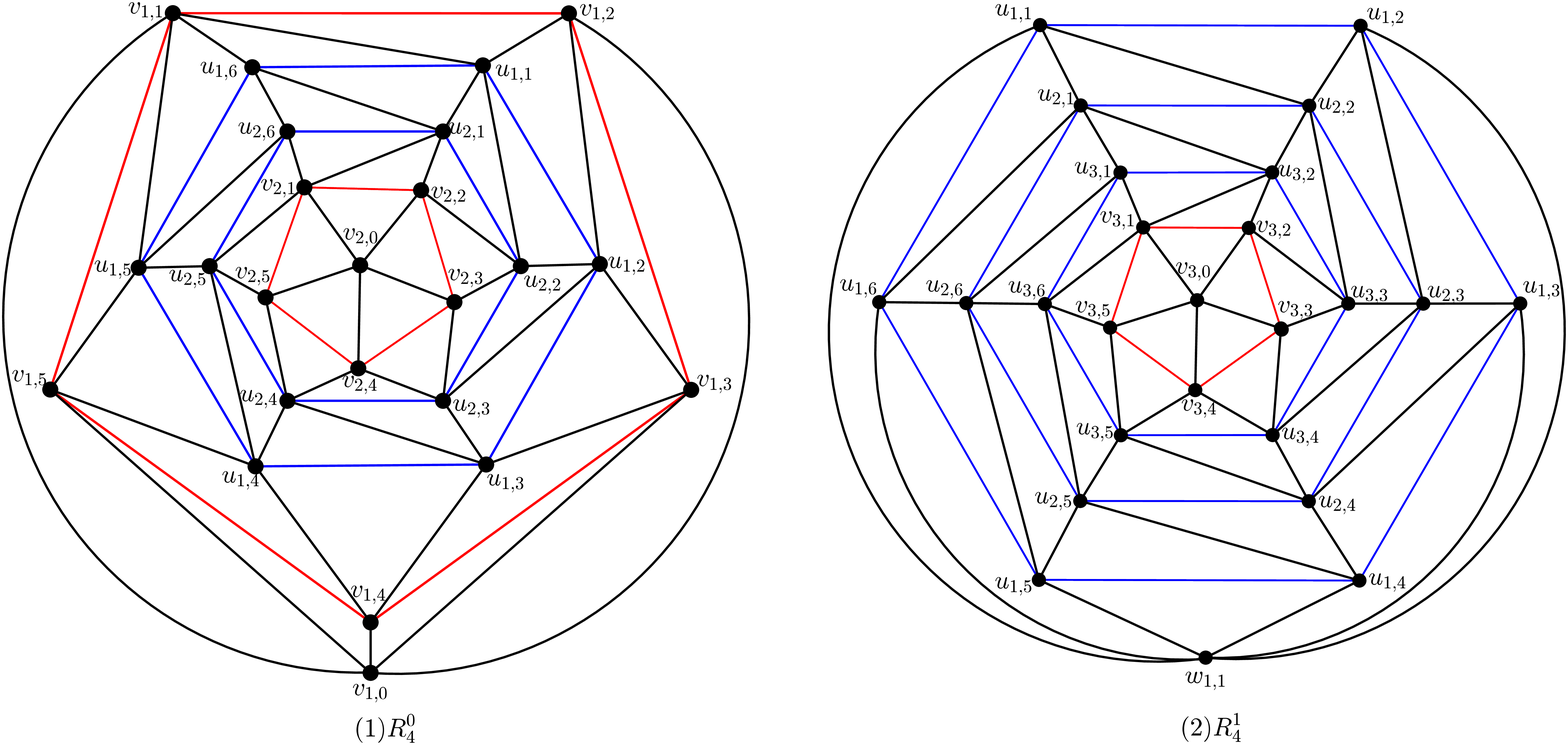}
\caption{The plane triangulations $R_k^0$ and $R_k^1$ when $k=4$, where $C^1,C^2,C^3$ are in  blue.}\label{R0}
\end{figure}

\begin{figure}[htbp]
\centering
\includegraphics*[scale=0.25]{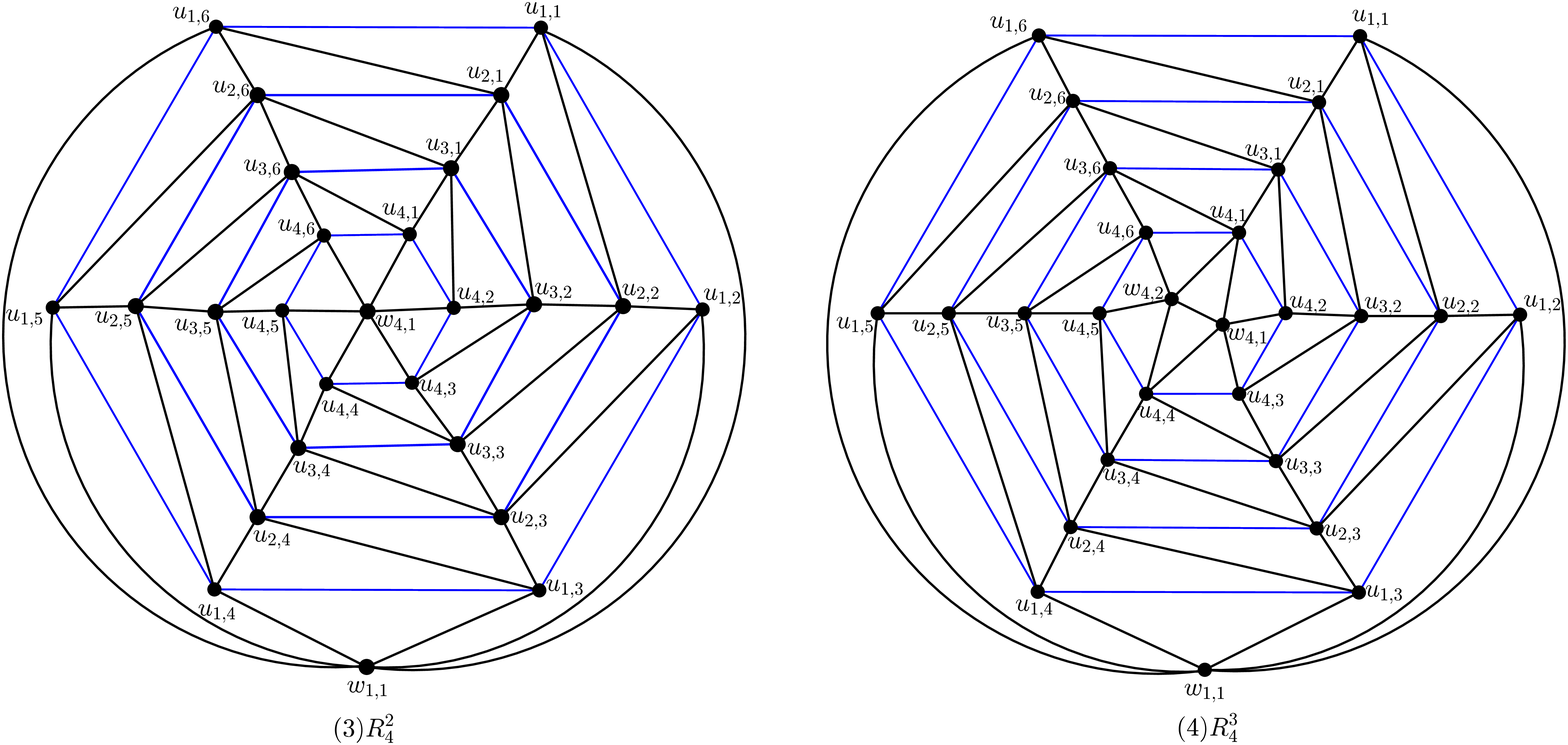}
\caption{The plane triangulations $R_k^2$ and $R_k^3$ when $k=4$, where $C^1,C^2,C^3,C^4$ are in  blue.}\label{R1}
\end{figure}

\begin{figure}[htbp]
\centering
\includegraphics*[scale=0.28]{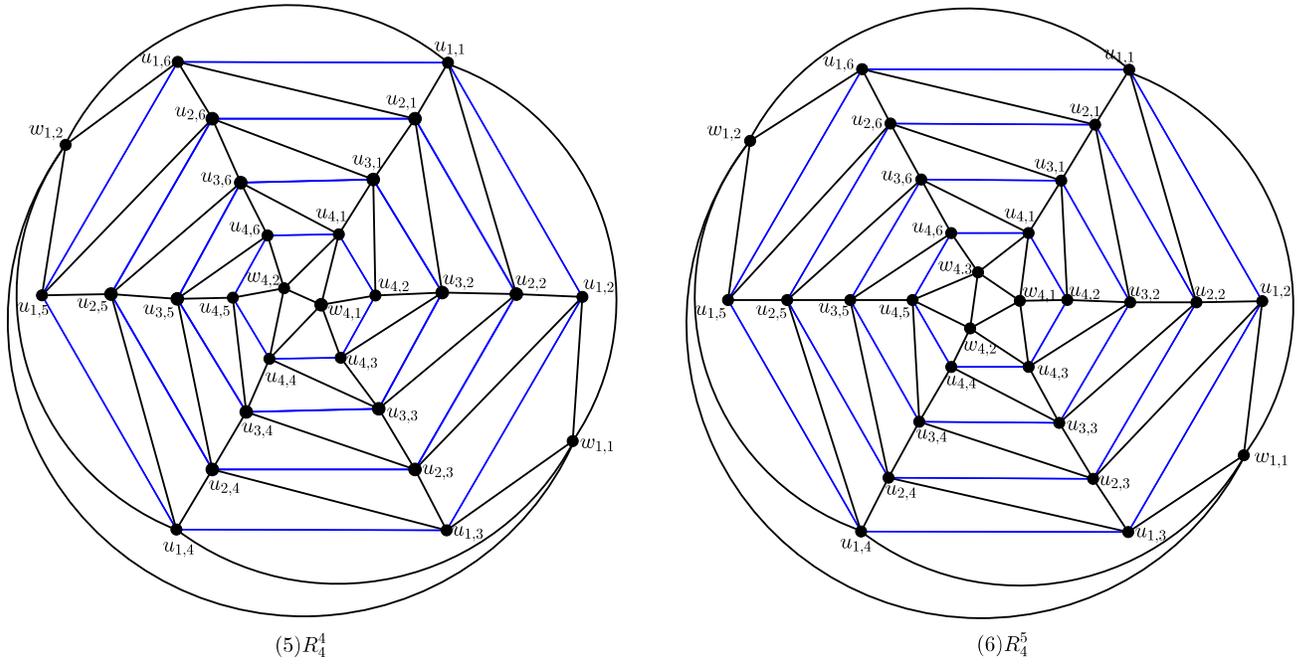}
\caption{The plane triangulations $R_k^4$ and $R_k^5$ when $k=4$, where $C^1,C^2,C^3,C^4$ are in  blue.}\label{R2}
\end{figure}

 Let  $n:=6k+r$ for some $r\in\{0, 1, \dots, 5\}$ and   integer $k\ge2$. When $r=0$,  let $R_k^0$ be the plane triangulation  on $n$ vertices obtained from $H_{k-2}\cup W_5^1\cup W_5^{k-2}$ by adding edges $u_{i,j}v_{i,j}$, $u_{i,j}v_{i,j+1}$ and $u_{i,6}v_{i,1}$ for all $i\in\{1,k-2\}$ and $j\in [5]$, where all arithmetic on the index $j+1$ here and henceforth is done modulo 5.
When $r=1$,  let $R_k^1$ be the plane triangulation on $n$ vertices obtained from $H_{k-1}\cup W_5^{k-1}\cup K_1^1$ by adding $w_{1,1}u_{1,m}$ for any $m\in[6]$ and adding edges $v_{k-1,j}u_{k-1,j}$, $v_{k-1,j}u_{k-1,j+1}$ and $u_{k-1,6}v_{k-1,1}$ for all $j\in [5]$, where all arithmetic on the index $j+1$ here is done modulo 5. When $r=2$, let $R_k^2$ be the plane triangulation on $n$ vertices obtained from $H_{k}\cup K_1^1\cup K_1^k$ by adding edges $u_{i,j}w_{i,1}$ for all $i\in\{1,k\}$ and $j\in[6]$. When $r=3$, let $R_k^3$ be the plane triangulation on $n$ vertices obtained from $H_{k}\cup K_1^1\cup K_2^k$ by adding edges $u_{1,i}w_{1,1}$ and $u_{k,j}w_{k,1}$ and $u_{k,s}w_{k,2}$ for all $i\in[6]$, $j\in[4]$ and $s\in\{1,4,5,6\}$. When $r=4$, let $R_k^4$ be the plane triangulation on $n$ vertices obtained from $H_{k}\cup K_2^1\cup K_2^k$ by adding edges $u_{i,j}w_{i,1}$ and $u_{i,s}w_{i,2}$ for all $i\in\{1,k\}$, $j\in[4]$ and $s\in\{1,4,5,6\}$. When $r=5$, let $R_k^5$ be the plane triangulation on $n$ vertices obtained from $H_{k}\cup K_2^1\cup K_3^k$ by first adding edges $u_{1,j}w_{1,1}$ and $u_{1,s}w_{1,2}$ for all $j\in[4]$ and $s\in\{1,4,5,6\}$, and then joining $w_{k,1}$ to  $u_{k, 1}, u_{k, 2},u_{k, 3}$; $w_{k,2}$ to   $u_{k, 3}, u_{k, 4},u_{k, 5}$;    $w_{k,3}$  to $u_{k, 1}, u_{k, 5},u_{k, 6}$.
For all  $r\in\{0,1,\ldots,5\}$, the graphs $R_k^r$ when $k=4$ are depicted in Figures~\ref{R0}-\ref{R2}. One can check that $R_k^{r}$  is $K_{2,t}$-free and so $ex_{_{\mathcal{P}}}(n, K_{2,t})=3n-6$ for  all $t\in\{3,4\}$ and $n\ge 12$.\medskip

We then consider the case $t=4$ and $6\le n\le 11$.
Let $Q_2$ be the plane graph defined in the proof of Theorem~\ref{main}. Let $Q_2'$ be the plane triangulation on $9$ vertices obtained from $Q_2\less v$ by adding edge $u_{2,1}u_{2,3}$; $Q_2''$ be the plane triangulation on $11$ vertices obtained from $Q_2$ by adding one vertex $w$ joining to $v,u_{2,1},u_{2,2}$. Then $Q_2',Q_2,Q_2''$ are $K_{2,4}$-free. Hence, $ex_{_{\mathcal{P}}}(n, K_{2,4})=3n-6$ for all $9\le n\le 11$. Let $G_1$ and $G_2$ be the unique plane graphs with degree sequence $5,4,4,3,3,3$ and $6,4,4,4,4,3,3$ respectively. Then the plane graphs $G_1,G_2,Q_2'\less u$ are $K_{2,4}$-free. Hence, $ex_{_{\mathcal{P}}}(n, K_{2,4})\ge3n-7$ for all $n\in\{6,7,8\}$.  Clearly, for $n\in\{6,7\}$, every plane triangulation on $n$ vertices is not $K_{2,4}$-free. Hence, $ex_{_{\mathcal{P}}}(n, K_{2,4})=3n-7$ for all $n\in\{6,7\}$. We next prove that $ex_{_{\mathcal{P}}}(n, K_{2,4})\le 3n-7$ when $n=8$. Suppose not. Let $G$ be a $K_{2,4}$-free plane triangulation on $n$ vertices. Let $v\in V(G)$ with $d_G(v)=\delta(G)$. Then $\delta(G)=4$, else $e(G)\ge{5n}/{2}>3n-6$ because $n=8$ when $\delta(G)\ge5$, or $e(G)=e(G\less v)+d_G(v)\le3(n-1)-7+3=3n-7$ when $\delta(G)\le3$. Then there exists $u\in V(G)\less N_G[v]$ such that $N_G(u)\cap N_G(v)\ge3$. But then $G$ contains a $K_{2,4}$ as a subgraph, a contradiction, as desired.  Hence, $ex_{_{\mathcal{P}}}(n, K_{2,4})=3n-7$ when $n=8$. \medskip

\begin{figure}[htbp]
\centering
\includegraphics*[scale=0.25]{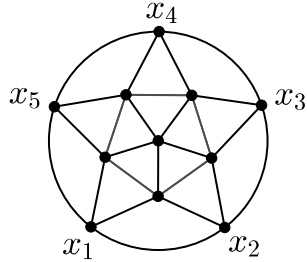}
\caption{The graph $J$.}\label{J}
\end{figure}

It remains to consider the case $t=3$ and $5\le n \le 11$.
Since $T_5^-$ has a copy of $K_{2,3}$ and $K_1+P_4$ is $K_{2,3}$-free, we have $ex_{_{\mathcal{P}}}(n, K_{2,3})=3n-8$ when $n=5$.
Let $O_7'$ be the near 4-regular plane graph obtained from $O_7$ by adding edges between vertices of degree $i$, where $i\in\{2,3\}$;
$O_8'$ be the 4-regular plane graph obtained from  $O_8$ by adding edges between vertices of degree at most 3. Let $J$ be the plane graph given in Figure~\ref{J}. Let $J'$ be the plane graph obtained from $J\less x_2$ by adding edge $x_1x_3$, $J''$ be the plane graph obtained from $J\less\{x_1,x_3\}$ by joining $x_2$ to $x_4,x_5$. Then the plane graphs $K_1+C_5,O_7',O_8',J'',J',J$ are $K_{2,3}$-free. Hence, $ex_{_{\mathcal{P}}}(n, K_{2,3})\ge3n-8$ for all $6\le n\le11$. We shall show that $ex_{_{\mathcal{P}}}(n, K_{2,3})\le3n-8$ for all $6\le n\le11$. Suppose this is not true. Let $G$ be a $K_{2,3}$-free plane graph on $n$ vertices with  $e(G)\ge 3n-7$, where $6\le n\le11$. We choose such a $G$ with $n$ minimum. Let $v\in V(G)$ with $d_G(v)=\delta(G)$. Then $\delta(G)\le4$, else $e(G)\ge{5n}/{2}>3n-6$ because $n\le 11$, a contradiction. Next, if $\delta(G)\le3$, then $e(G\less v)\le 3(n-1)-8$ by the minimality of $n$ and the fact that $ex_{_{\mathcal{P}}}(n, K_{2,3})\le3n-8$ when $n=5$. Thus, $e(G)=e(G\less v)+d_G(v)\le 3(n-1)-8+3=3n-8$, a contradiction. This proves that $\delta(G)=4$. We see  $G$ is not plane triangulation because $G$ is $K_{2,3}$-free. So $G$ has at least three vertices of degree four, else $e(G)\ge\lceil{(5n-2)}/{2}\rceil\ge3n-6$ because $n\le11$, a contradiction. This implies that $G$ contains $K_1+C_4$ as a subgraph and so $G$ is not $K_{2,3}$-free, a contradiction. Therefore, $ex_{_{\mathcal{P}}}(n, K_{2,3})=3n-8$ for all $6\le n\le11$.\medskip

This completes the proof of Theorem~\ref{bipartite}.

\section{Concluding remarks}

Theorem~\ref{tC3} completely determines the planar Tur\'an number of $2C_3$.  It seems non-trivial to determine $ex_{_\mathcal{P}}(n,2C_k)$ for all $k\ge 4$ even with the lower bound given in Lemma~\ref{2Ck}.
However,  the lower bound in Lemma~\ref{2Ck}  seems tight when $k=4$.
\begin{conj}\label{C4}
For all $n\ge 8$,  let  $r$ be the remainder of $n-3$ when divided by $2$. Then
\[ex_{_\mathcal{P}}(n,2C_4)=\frac{5n}{2}-\frac{5+r}{2}.\]
 \end{conj}
It   seems that  Proposition~\ref{p:C3+} can be extended to all $k\ge4$, that is,    $ex_{_{\mathcal{P}}}(n,  C_k^+)=ex_{_{\mathcal{P}}}(n,C_k)$ for all $n\ge k+1\ge5$.  We are unable to present a proof here, even with the assumed truth of Conjecture~\ref{weak}. \medskip

We end the paper with a lemma which improves  the lower bound given in Lemma~\ref{2Ck}   for all  $k\ge 7$. It would be interesting to know whether the bounds in Lemma~\ref{2Cknew} below are desired  tight upper bounds for $ex_{_\mathcal{P}}(n,2C_k)$ for all  $k\ge 7$ and $n$ sufficiently large.

 \begin{lem}\label{2Cknew}
Let $n$ and $k$ be positive integers with  $n\ge 2k\ge 14$. Let $\varepsilon_1$   and $\varepsilon_2$ be the remainder of $n-(2k-1)$ when divided by $k-4+\frac{k-1}{2}$ ($k$ is odd)  and   $k-6+\frac{k}{2}$ ($k$ is even), respectively. 

 \begin{enumerate}[(a)]
\item If  $k$  is  odd,  then  $ex_{_\mathcal{P}}(n,2C_k)=3n-6$ for all  $n\le 3k-4$, and
\[ex_{_\mathcal{P}}(n,2C_k)\ge\left(3-\frac{1}{k-4+\lfloor  k/2\rfloor}\right)n+\frac{5+\varepsilon_1}{k-4+\lfloor  k/2\rfloor}-\frac{17}{3}+\max\{1-\varepsilon_1, 0\}\,  \text{ for all }  n\ge 3k-3.\]
\item If $k$  is   even,   then   $ex_{_\mathcal{P}}(n,2C_k)=3n-6$ for all $  n\le 3k-7$, and  
\[ex_{_\mathcal{P}}(n,2C_k)\ge\left(3-\frac{1}{k-6+ k/2}\right)n+\frac{7+\varepsilon_2}{k-6+ k/2}-\frac{17}{3}+\max\{1-\varepsilon_2, 0\} \,  \text{ for all }  n\ge 3k-6.\]
 \end{enumerate}
\end{lem}
\pf  Let $n,k,\varepsilon_1, \varepsilon_2$ be given as in the statement. Throughout the proof, let $\mathcal{T}_m:=K_2+P_{m-2}$ be the plane triangulation on $m\ge2$ vertices, with $x$ and $y$ on the outer face of $\mathcal{T}_m$, where $x$ and $y$ are the two adjacent vertices of degree $m-1$ in $\mathcal{T}_m$. Note that we allow $m$ to be $2$ here for a simpler proof later on; $\mathcal{T}_m$ has exactly $2m-4$ $3$-faces. For each integer $s$ satisfying $m\le  s\le 2m-4$, let $\mathcal{T}_s^m$ denote a plane triangulation on $s$ vertices obtained from $\mathcal{T}_m$ by adding $s-m$ new vertices:  each to a $3$-face $F$ of $\mathcal{T}_m$ and then joining it to all vertices on the boundary of $F$.  \medskip

To prove ($a$),  let $k:=2p+1$, where $p\ge3$ is an integer. We first show that $ex_{_\mathcal{P}}(n,2C_k)=3n-6$ for all $n\le 3k-4$.  Note that $\mathcal{T}_{2p+1}$ has exactly $2(2p+1)-4=4p-2$ 3-faces. For each $n\le 3k-4=6p-2=(2p+1)+(4p-2)$,    $\mathcal{T}^{2p+1}_n$ is $2C_k$-free because each $C_k$ in $\mathcal{T}^{2p+1}_n$ must contain at least $p+1$ vertices of $\mathcal{T}_{2p+1}$. Hence,  $ex_{_\mathcal{P}}(n,2C_k)=3n-6$ for all $n\le 3k-4$. We next consider the case $n\ge3k-3$.
Let $t\ge0$ be an integer satisfying
$$t(3p-3)+\varepsilon_1=n-(2k-1).$$
Let  $H_i:=\mathcal{T}^{p+1}_{3p-1}$ for  each  $i\in[t]$ when $t\ge1$;    $H_{t+1}:=\mathcal{T}_{\varepsilon_1+2}$  when $ \varepsilon_1+2\le k-1$ and   the plane triangulation   $\mathcal{T}^{p+1}_{\varepsilon_1+2}$  when $\varepsilon_1+2\ge k$; and  $H_{t+2}:=\mathcal{T}_{2k-1}$.     For each $j\in[t+1]$, it is easy to check that   each cycle of $H_j$ has length at most $k$; each cycle of length exactly $k$ must contain the vertices $x$ and $y$; $H_j$ is  $2C_k$-free.  Note that $H_{t+2}$ is $2C_k$-free,  and each cycle of length exactly $k$ in $H_{t+2}$ must contain  either  $x$ or  $y$ (or possibly both).   Finally, let $G$ be the  planar graph obtained from disjoint copies of $H_1, H_2, \ldots, H_{t+2}$ by pasting along the subgraph $K_2$ induced by $x$ and $y$.  It follows that $G$ is $2C_k$-free and  $|G|=t(3p-3)+(2k-1)+\varepsilon_1=n$.  Note that $|H_{t+2}|\ge 2$ with equality when  $\varepsilon_1=0$.   Therefore,   
\begin{align*}
ex_{_\mathcal{P}}(n,2C_k)\ge e(G)&=\sum_{i=1}^t(e(H_i)-1)+(e(H_{t+1})-1)+(e(H_{t+2})-1)+1\\
&=t[3(3p-1)-7]+ \max\{3(\varepsilon_1+2)-7, 0\}+(3(2k-1)-7)+ 1\\
&=3n-t-7+\max\{1-\varepsilon_1, 0\} \\
&=3n-\frac{n- (2k-1)-\varepsilon_1}{3p-3}-7+\max\{1-\varepsilon_1, 0\}\\
&=\left(3-\frac{1}{3p-3}\right)n+\frac{2k-1+\varepsilon_1}{3p-3}-7+\max\{1-\varepsilon_1, 0\}\\
&=\left(3-\frac{1}{k-4+\frac{k-1}2}\right)n+\frac{2k-1+\varepsilon_1}{k-4+\frac{k-1}2}-7+\max\{1-\varepsilon_1, 0\}\\
&=\left(3-\frac{1}{k-4+\lfloor  k/2\rfloor}\right)n+\frac{5+\varepsilon_1}{k-4+\lfloor  k/2\rfloor}-\frac{17}{3}+\max\{1-\varepsilon_1, 0\}.\\
 \end{align*}

 It remains to prove ($b$).    Let $k: =2p$, where $p\ge4$ is an integer. Similar to the proof of ($a$), we see that $\mathcal{T}^{2p-1}_n$ is $2C_k$-free  for each $n\le 3k-7$, because each $C_k$ in $\mathcal{T}^{2p-1}_n$ must contain at least $p$ vertices of $\mathcal{T}_{2p-1}$. Hence,  $ex_{_\mathcal{P}}(n,2C_k)=3n-6$ for all $n\le 3k-7$. We next consider the case $n\ge3k-6$. 
Let $t\ge0$ be an integer satisfying
$$t(3p-6)+\varepsilon_2=n-(2k-1).$$
Let  $H_i:=\mathcal{T}^{p}_{3p-4}$ for  each  $i\in[t]$ when $t\ge1$;   $H_{t+1}:=\mathcal{T}_{\varepsilon_2+2}$  when $ \varepsilon_2+2\le k-1$ and  the plane triangulation   $\mathcal{T}^{p}_{\varepsilon_2+2}$  when $\varepsilon_2+2\ge k$; and $H_{t+2}:=\mathcal{T}_{2k-1}$.    For each $j\in[t+1]$, it is easy to check that   each cycle of $H_j$ has length at most $k$; each cycle of length exactly $k$ must contain the vertices $x$ and $y$; $H_j$ is  $2C_k$-free.  Note that $H_{t+2}$ is $2C_k$-free,  and each cycle of length exactly $k$ in $H_{t+2}$ must contain  either  $x$ or  $y$ (or possibly both).   Finally, let $G$ be the  planar graph obtained from disjoint copies of $H_1, H_2, \ldots, H_{t+2}$ by  pasting along the subgraph $K_2$ induced by $x$ and $y$.  It follows that $G$ is $2C_k$-free and  $|G|=t(3p-6)+(2k-1)+\varepsilon_2=n$.  Note that $|H_{t+2}|\ge 2$ with equality when  $\varepsilon_2=0$.   Therefore,   
\begin{align*}
ex_{_\mathcal{P}}(n,2C_k)\ge e(G)&=\sum_{i=1}^t(e(H_i)-1)+(e(H_{t+1})-1)+(e(H_{t+2})-1)+1\\
&=t[3(3p-4)-7]+\max\{3(\varepsilon_2+2)-7, 0\}+(3(2k-1)-7)+1\\
&=3n-t-7+\max\{1-\varepsilon_2, 0\} \\
&=3n-\frac{n- (2k-1)-\varepsilon_2}{3p-6}-7+\max\{1-\varepsilon_2, 0\}\\
&=\left(3-\frac{1}{3p-6}\right)n+\frac{2k-1+\varepsilon_2}{3p-6}-7+\max\{1-\varepsilon_2, 0\}\\
&=\left(3-\frac{1}{k-6+ k/2}\right)n+\frac{2k-1+\varepsilon_2}{k-6+ k/2}-7+\max\{1-\varepsilon_2, 0\}\\
&=\left(3-\frac{1}{k-6+ k/2}\right)n+\frac{7+\varepsilon_2}{k-6+ k/2}-\frac{17}{3}+\max\{1-\varepsilon_2, 0\}, 
 \end{align*}
 as desired. This completes the proof of Lemma~\ref{2Cknew}.\qed

\medskip

%

\end{document}